\documentclass[12pt]{article}
\usepackage{latexsym,amsthm,amsmath}
\def\dj{d\kern-.30em\raise1.25ex\vbox{\hrule width .3em height .03em}}
\def\Dj{D\rlap{\kern-.70em\raise0.75ex
\vbox{\hrule width .3em height .03em}}}

\def\rar{\rightarrow}
\def\bk{I \kern-.25em k}
\def\nd{{\mathrm{End}}\,}
\def\hom{{\mathrm{Hom}}\,}
\def\mod{{\mathrm{Mod}}\,}
\def\id{\mathrm{id}}
\def\dim{\mathrm{dim}}
\def\de{\delta}
\def\ep{\varepsilon}
\def\al{\alpha}

\def\lam{\lambda}
\def\Lam{\Lambda}
\def\pt{\partial}
\def\yd{{\mathcal Y}{\mathcal D}}
\def\co{\triangle}
\def\cl#1{\bigtriangleup_{#1}}
\def\ccr#1{_{#1}\!\bigtriangleup}

\def\op{{\mathrm op}}
\def\cop{{\mathrm cop}}
\def\ot{\otimes}
\def\c{\circ}
\def\rfr#1#2{#1_{#2}}
\def\lfr#1#2{_{#1}\! #2}
\def\rd#1{#1^\dagger}
\def\ld#1{^\dagger\! #1}
\newtheorem{defn1}{Definition}[section]
\newtheorem{rem1}[defn1]{Remark}
\newtheorem{prop1}[defn1]{Proposition}
\newtheorem{cor1}[defn1]{Corollary}
\newtheorem{prop2}[defn1]{Proposition}
\newtheorem{rem2}[defn1]{Remark}
\newtheorem{lem1}{Lemma}[section]
\newtheorem{rem3}[lem1]{Remark}
\newtheorem{rem4}[lem1]{Remark}
\newtheorem{prop3}[lem1]{Proposition}
\newtheorem{rem5}[lem1]{Remark}
\newtheorem{rem6}[lem1]{Remark}
\newtheorem{lem2}{Lemma}[section]
\newtheorem{lem3}[lem2]{Lemma}
\newtheorem{thm}[lem2]{Main Theorem}
\newtheorem{cor2}[lem2]{Corollary}
\newtheorem{prop4}[lem2]{Proposition}
\begin{document}
\title{Hopf Modules and Their Duals\thanks{Supported by
Polish KBN grant \# 2 P03B 109 15; by UNAM, DGAPA, Programa de
Apoyo a Proyectos de Investigacion e Innovacion Tecnologica,
M\'exico proyecto \# IN-109599; and by el Consejo Nacional de
Ciencia y Tecnolog\'{\i}a CONACyT de M\'exico, proyecto \#
27670E.}}
\date{to appear in Int. J. Theoret. Phys. v.39 no10 (2000)}
\author{Andrzej Zdzis{\l}aw Borowiec\thanks{University of Wroc{\l}aw,
Institute of Theoretical Physics, plac Maksa Borna 9, PL 50-204
Wroc{\l}aw, Poland, borow@ift.uni.wroc.pl}\\
Guillermo Arnulfo Va\'zquez Couti\~no\thanks{Universidad
Aut\'onoma Metropolitana-Iztapalapa, Apartado Postal 55-534,
M\'exico D.F., C.P. 09340, gavc@xanum.uam.mx}}
\maketitle
\begin{abstract}
Free Hopf modules and bimodules over a bialgebra are studied with
some details. In particular, we investigate a duality in the
category of bimodules in this context. This gives the
correspondence between Woronowicz's quantum Lie algebra and
algebraic vector fields.\smallskip\\ \noindent\textbf{1999 Physics
and Astronomy Classification Scheme:} 02.10.Tq, 02.20.-a,
02.40.k\\ \noindent\textbf{2000 Mathematics Subject
Classification}: Primary 16W30, 17B37, Secondary 81R50, 16S40.\\
\end{abstract}
\section{INTRODUCTION}
In this note, we are interested in comparing the Woronowicz
quantum Lie algebra construction for a quantum group bicovariant
differential calculus \cite{Wor} with the algebraic vector fields
for the first order differential calculus over an arbitrary unital
associative algebra \cite{Bor1,Bor2}. Our result is that dualizing
a bicovariant bimodule of one-forms in the category of bimodules
over a Hopf algebra with bijective antipode, one obtains bimodule
of algebraic vector fields. Like in the Lie algebra case, the
quantum Lie algebra consists of the left or right invariant vector
fields. This bimodule acts as a \textit{Cartan pair} \cite{Bor1}
on the algebra itself, and it is bicovariant over co-opposite Hopf
algebra. In particular, one questions statement formulated in
\cite{Asc,AS} that it is bicovariant over the same Hopf algebra.

Since Woronowicz celebrated paper \cite{Wor}, bicovariant
differential calculi have become a subject of a huge number of
investigations, see e.g. \cite{Ber,KS,Maj,Ozi,Dur}. General
construction of vector fields for bicovariant differential calculi
on Hopf algebras have been previously discussed by Aschieri \&
Schupp \cite{AS,Asc}, Pflaum \& Schauenburg \cite{PS} and
Schauenburg \cite{Sch2}. In fact, our construction generalize that
of \cite{AS}.

The problem of generalization of the Lie module to quantum and
braided category is still open. Among several propositions we like
to mention the Pareigis approach to Lie algebras in the category
of Yetter-Drinfeld modules \cite{Par}, and several approaches
based on variety of braided identities generalizing Jacobi
identity \cite{OPR,OR,Bau}.

In the sequel $\bk$ is a field. We shall work in the category of
$\bk$-vector spaces, all maps are $\bk$ linear maps and the tensor
product is over $\bk.$ Given $\bk$-spaces $U$ and $W$, $\hom(U,W)$
denote the $\bk$-space of all $\bk$-linear maps from $U$ to $W$.
Denote by $\tilde U\doteq\hom(U,\bk)$ the $\bk$-linear dual of
$U$. For $\Phi\in\hom(U,W)$ we denote by $\tilde\Phi\in\hom(\tilde
W, \tilde U)$ its linear transpose. Dealing with
finite-dimensional vector spaces we shall use the covariant index
notation together with the Einstein summation convention over
repeated up contravariant and down covariant indices. If
$\{e_k\}_{k=1}^{{\mathrm dim}V}$ denotes a basis in a
finite-dimensional $\bk$-space $V,$ then $v=v^i e_i\in V.$

An algebra means associative unital $\bk$-algebra and a coalgebra
means coassociative counital $\bk$-coalgebra. If $A \equiv(A,m,1)$
is an algebra then $A^\op$ denotes algebra with the opposite
multiplication: $a\cdot_\op b=b\,a$ Let $C\equiv (C,\co,\ep$) be
any coalgebra with comultiplication $\co$ and counit $\ep$. The
Sweedler \cite{Swe} shorthand notation is $\co (a)= a_{(1)}\ot
a_{(2)}$. For a left (right) comultiplication on $V$ we shall
write $\co_V (v)=v_{(-1)}\ot v_{(0)}$ ($_V\co (v)=v_{(0)}\ot
v_{(1)}$ respectively). By $C^\cop$ we mean a opposite coalgebra
structure with $\co^\cop (a)= a_{(2)}\ot a_{(1)}$. For a given
bialgebra $B$, one can form  new bialgebras by taking the opposite
of either the algebra or/and coalgebra structure, e.g. $B^{\op\
\cop}$ has both opposite structures.

Various modules and comodules over a bialgebra are our main
objects of investigation. The most substantial results are
obtained for the case of bialgebras with the bijective antipode
(quantum groups).

\section{NOTATION AND PRELIMINARIES} Let $A$ be an algebra.
Assume that a finite dimensional $\bk$-space $V$ is a left
$A$-module, or equivalently, it is a carrier space for
representation $\lam$ of $A$. This means that the left action
$m_V:A\ot V\rar V$ can be written in terms of a unital
$\bk$-algebra homomorphism $\lam:A\rar\nd V$ which in turn, by the
use of an arbitrary basis $\{e_k\}_{k=1}^{\dim V}$ of $V$, can be
rewritten in matrix form $\lam^i_k(a)e_i\doteq m_V(a\ot e_k)\doteq
\lam(a)e_k,$ $$\forall\,a,b\in A,\quad \lam^i_k(1)=\de^i_k,\quad
\lam^i_k(ab)= \lam^i_m(a) \lam^m_k(b)\ . \eqno (1) $$ The same
matrix representation uniquely induces the transpose right
multiplication $_{\tilde V}\!m\doteq\widetilde{m_V}: \tilde V\ot
A\rar\tilde V$ on the dual vector space $\tilde V$: $$
\widetilde{m_V}(e^k\ot a)\doteq \tilde\lam(a)e^k\equiv
\lam^k_m(a)e^m\ ,\eqno (2)$$ where $e^k$-s are elements of the
dual basis in $\tilde V$. This defines an anti-representation
$\tilde\lam:A\rar\nd\tilde V$ given by the transpose matrices
$\tilde\lam(a)$: $\tilde\lam(ab)=\tilde\lam(b)\tilde\lam(a)$.
Alternatively, one can say that  $\tilde\lam: A^{\op}\rar
\nd\tilde V$ is a representation of the opposite algebra $A^{\op
}$, i.e. it defines a left $A^{\op}$ module structure on $\tilde
V$.

For an algebra morphism $T:A^\prime\rar A$ and a left $A^\prime$-
module action $m^\prime_L: A^\prime\ot V\rar V$ one defines its
{\it pull-back} as a left $A$-module action $m_L\doteq T^\star
(m^\prime_L): A\ot V\rar V$ by $$T^\star (m^\prime_L)\doteq
m^\prime_L\c (T\ot\id),\eqno(3)$$ $$a\cdot_Tv\doteq
T(a).v\quad\hbox{or}\quad\lam^k_i\doteq \lam^{\prime k}_i\c T.$$
When $T$ is an anti-homomorphism (i.e. a homomorphism from $A$
into $A^{\prime\ \op}$ ) then the pull-back $T^\star (m^\prime_L)$
is a right $A$-module action.

 Let $C$ be a coalgebra. A left $C$-comodule structure on
$V$ is determined by {\it matrix-like} elements $L^i_k\in C$,
$i,k=1,\ldots ,\dim V$, $$ \co (L^i_k)= L^m_k\ot L^i_m,\ \ \ \ep
(L^i_k)=\de^i_k .\eqno (4) $$ These define the left coaction or
corepresentation $\cl V:V\rar C\ot V$: $$ \cl
V(e_k)=(e_k)_{(-1)}\ot(e_k)_{(0)}\doteq L^m_k\ot e_m . \eqno (5)$$
The same matrix elements $L^i_k\in C$ induce the transpose right
comultiplication $\ccr{\tilde V}\doteq \widetilde{\cl V}:\tilde
V\rar\tilde V\ot C$ $$\widetilde{\cl
V}(e^k)=(e^k)_{(0)}\ot(e^k)_{(1)}\doteq e^m\ot L^k_m \eqno (6)$$
on the dual vector space $\tilde V$ (cf. \cite{BV}).
Alternatively, one can say that $\widetilde{\cl V}$ defines a left
coaction of the co-opposite coalgebra $C^{\mathrm cop}$ on $\tilde
V$.

For a coalgebra morphism $T:C\rar C^\prime$ and a left
$C$-comodule coaction $\cl L: V\rar C\ot V$ one defines its {\it
push-forward} $\cl L^\prime\equiv T_\star (\cl L): V\rar
C^\prime\ot V$ as a left $C^\prime$-comodule  coaction such that
$$T_\star(\cl L)\doteq (T\ot\id)\c\cl V\ \ \ \ \hbox{or}\ \ \ \
L^{\prime k} _i=T(L^k_i) .\eqno (7)$$ If $T$ is an anti-coalgebra
map then its push-forward $T_\star (\cl L)$ is a right
$C^\prime$-coaction on $V$.

\section{YETTER-DRINFELD MODULES}
Our basic references on Yetter-Drinfeld ($\yd$) modules are
\cite{Mon,RT,Sch1}. $\yd$ modules are also known under the name of
Yang-Baxter or crossed modules. Let $V$ be some finite dimensional
$\bk$-space.

\begin{defn1} For bialgebra B, a left-left {\it Yetter-Drinfeld}
module is a $\bk$-space $V\equiv (V, m_V, \cl V)$ which is both a
left $B$-module and a left $B$-comodule and satisfies the
compatibility condition $$\forall\,a\in B\;\text{and}\;v\in
V,\quad a_{(1)}v_{(-1)}\ot a_{(2)}v_{(0)} =
(a_{(1)}.v)_{(-1)}a_{(2)}\ot(a_{(1)}.v)_{(0)},$$ $$\text{or
equivalently} \ \ \ \ \
a_{(1)}L_k^m\lam^i_m(a_{(2)})=L^i_ma_{(2)}\lam^m_k(a_{(1)})\
,\eqno (8) $$ 
\end{defn1}

We denote by $_B^B\yd$ the category of left-left $\yd$ modules
over $B$. Similarly, the right-right $\yd$ module condition is
$$v_{(0)} a_{(1)}\ot v_{(1)}a_{(2)} = (v.a_{(2)})_{(0)}\ot
a_{(1)}(v.(a_{(2)}))_{(1)},$$ $$\text{or}\ \ \ \ \
a_{(1)}R^k_m\rho^m_i(a_{(2)})=R^m_ia_{(2)}\rho^k_m(a_{(1)})\
,\eqno (9) $$ where $_V\!m(e_k\ot a)\doteq \rho^m_k(a)e_m$ and
$_V\!\co (e_k)\doteq e_m\ot R^m_k$ denotes the right
multiplication and comultiplication in $V$. Consequently,
$\yd^B_B$ denotes the category of right-right $\yd$ $B$-modules.
Likewise, one can introduce categories of left-right $_B\yd^B$ and
right-left $^B\yd_B$ Yetter-Drinfeld modules over $B$.

\begin{rem1} One recognizes that a left-left $\yd$ module over a
bialgebra $B$ becomes automatically a right-right $\yd$ module
over the bialgebra $B^{\op\ \cop}$. More generally, the following
categories $$_B^B\yd\equiv\ _{B^\cop}\!\yd^{B^\cop}\equiv\
^{B^\op}\!\yd_{B^\op}\equiv\yd^{B^{\op\ \cop }}_{B^{\op\ \cop }}
\eqno (10)$$ can be identified in the formal sense, i.e.: if a
triple $(V, m_L, \cl L)$ is an element of the first category then
it becomes automatically (after suitable re-interpretations) an
element of the remaining categories. E.g., denoting by $\cl L^\cop
(e_k) =e_m\ot L^m_k$ a canonical right $B^\cop$-comodule structure
on $V$ associated with $\cl L$, one sees that $(V, m_L, \cl
L^\cop)\in\, _{B^\cop}\!\yd^{B^\cop}$.
\end{rem1}

Nontrivial category equivalences have been found for the special
case, when $B$ is a Hopf algebra with bijective antipode.

\begin{prop1}[Radford \& Towber \cite{RT}, p. 265] Suppose that $B$ is a
bialgebra with bijective antipode $S$. Then
\begin{description}
\item[{\em i) (Woronowicz \cite{Wor})}] $(V, m_L, \cl L)\mapsto (V, (S^{-1})^\star
(m_L), S_\star (\cl L))$ describes categorical isomorphisms
$$^B_B\yd\cong\ \yd^B_B\quad\hbox{and}\quad _B\yd^B\cong \
^B\yd_B. \eqno (11)$$

\item[{\em ii)}] $(V, m_L, \cl R)\mapsto (V, m_L, S_\star (\cl R))$ describes
categorical isomorphisms $$_B\yd^B\cong\ _B^B\yd .\eqno (12)$$
\end{description}\end{prop1}

\begin{cor1} Combining (10) and (12) one gets the categories
equivalence $$^B_B\yd\ni(V,m_L, \cl L)\longmapsto (V, m_L,
(S^{-1})_\star (\cl L^\cop))\in\, ^{B^\cop}_{B^\cop}\yd \ ;$$
$$\yd^B_B\ni(V,m_R, \cl R)\longmapsto (V, m_R, (S^{-1})_\star (\cl
R^\cop))\in\, \yd^{B^\cop}_{B^\cop}. \eqno (13)$$
\end{cor1} These follow from the fact that $S^{-1}$ is an
antipode of $B^\cop$.

\begin{prop2} Let $(V, m_L, \cl L)$ be  a left module and left comodule over a
bialgebra $B$, with $\dim V<\infty$. Then $(V, m_L, \cl L)\in\
^B_B\yd$ if and only if $(\tilde V, \widetilde{m_L},
\widetilde{\cl L})\in\yd^B_B$.
\end{prop2}
\begin{proof} Substituting $L^i_k=R^k_i$ and $\lam^i_k=\rho^k_i$ into
equation (8) one gets (9). A right-handed version of the
Proposition also holds.\end{proof}  

\begin{rem2} A $\yd$ structure on $V$
 generates the quantum Yang-Baxter operator ${\mathcal R}\in \nd{(V\ot
 V)}$. For example, if $(V,\, _Vm,\, \ccr V)\in\, \yd^B_B$ one gets
 $${\mathcal R}(e_i\ot e_k)=\rho^j_i(R^m_k)\,e_m\ot e_j.$$
 \end{rem2}

\section{FREE COVARIANT BIMODULES} A left (resp. right) free $A$-module
$M$ can be represented as  $A\!\ot\! V$ (resp. $V\!\ot\! A$),
where $V$ denotes a linear space spanned by free generators
$\xi_1,\ldots ,\xi_n$, $n=\dim V$ and the module structure is
realized by the left (resp. right) multiplication in $A$.
Following Sweedler \cite{Swe} we shall use the notation $\rfr
VA\doteq V\ot A$ and $\lfr AV\doteq A\ot V$ for the right and the
left free $A$-modules generated by a vector space $V$. In the
sequel   we shall restrict ourselves exclusively to the case when
$V$ is a finite dimensional vector space. We do not assume an {\it
invariant basis property} for algebra $A$. This means that the
number of free generators is not necessarily a characteristic
number for a given free (left or right) module. In other  words,
one can have a left $A$-module isomorphism $A\ot V\cong A\ot W$
with $\dim V\neq \dim W$.

The unit $1\equiv 1_A$ of $A$ enables us to define a canonical
inclusion $V\ni v\mapsto v\ot 1\in\rfr VA$ and $$\forall\,x\in\rfr
VA,\quad x=e_i\ot x^i,\eqno (14)$$ where components $x^i\in A$ are
uniquely determined with respect to a given basis $\{e_i\}$. Any
basis $\{e_i\}$ in $V$ determines a set of free generators
$\{\xi_i=e_i\ot 1_A\}$ in the module $\rfr VA$.

\begin{lem1}[see e.g. \cite{Sch1,BKO}]
Let $V$ be a finite-di\-m\-en\-sio\-nal $\bk$-space. The following
are equivalent\begin{description}
\item[i)] A left $A$-module structure on a free right module $\rfr VA$
such that it becomes an $A$-bimodule.
\item[ii)] A unital $\bk$-algebra map (so called {\it commutation rule})
$\Lam:A\rar A\ot\nd V.$
\item[iii)] A $\bk$-linear map (so called {\it
twist}) $\hat\Lam:\,\lfr AV\rar\,\rfr VA$ such that
$$\hat\Lam(1\ot v)=v\ot 1,$$ $$\hat\Lam\c(m\ot\id_V)=(\id_V\ot
m)\c(\hat\Lam\ot\id_A)\c(\id_A\ot\hat\Lam)\ . $$\end{description}
\end{lem1}
\begin{proof} By uniqueness of the decomposition (14) one can set
$$\Lam(a)e_k\doteq a.(e_k\ot 1)\doteq \hat\Lam(a\ot e_k)\doteq
e_i\ot\Lam^i_k(a) \eqno (15)$$ in an arbitrary basis $\{e_i\}$ of
$V$. Properties of the $\bk$-algebra map $$\Lam^i_k(1)=\de^i_k,\ \
\ \Lam^i_k(ab)=\Lam^i_m(a)\Lam^m_k(b)\eqno (16)$$ as well as
$iii)$ are to be verified.\end{proof}

A right-handed version of the Lemma above: a right commutation
rule $\Phi:A\rar A\ot\nd V$  gives rise to a right multiplication
on $\lfr AV,$ $$(1\ot e_k).a\doteq\Phi_k^i(a)\ot e_i,\ \
\Phi^i_k(ab)=
  \Phi^m_k(a)\Phi^i_m(b),\ \ \Phi^i_k(1)=\de^i_k .\eqno (17)$$
This implies that $\Phi$ is an algebra map from $A^\op$ into
$A^\op\ot\nd V$.

If $C$ is a coalgebra then $\lfr CV$ is a left free comodule with
a comultiplication determined by that in $C$, i.e. $\cl{\lfr
CV}=\co\ot\id_V$. The counit $\ep$ in $C$ enables us to define a
projection  map $\ep_V:\ \lfr CV\rar V$ by $\ep_V(x^i\ot
e_i)\doteq \ep(x^i)e_i,$ $$\ep_V(a.x)=\ep(a)\ep_V(x) ,\ \ \
(\id\ot\ep_V)\c\cl{\lfr CV}=\id .\eqno (18)$$

Let $B$ be a bialgebra. In this case $\lfr BV$ is a (free) left
module and a (free) left comodule with multiplication and
comultiplication satisfying the following compatibility condition
$$\forall\,a\in B\;\text{and}\;x\in\,\lfr BV,\quad
 \cl{\lfr BV}(a.x)=\co(a)\cl{\lfr BV}(x),\eqno (19)$$
$$(a.x)_{(-1)}\ot (a.x)_{(0)} =a_{(1)}.x_{(-1)}\ot
a_{(2)}x_{(0)}.$$ This condition  differs from the $\yd$
conditions and defines a left {\it Hopf} $B$-module structure on a
left free module $\lfr BV$. Similarly, $\rfr VB$ becomes
automatically, a right free Hopf $B$-module.

\begin{rem3} $B$ and $B^{\cop}$ have the same algebra
structure. Therefore we can treat $\lfr BV$ as a left free Hopf
$B^{\cop}$-module with coaction $\co^{\cop}_{\lfr
BV}=\co^{\cop}\ot\id_V$. A $\bk$-space $V$ generates a free left
(right) either $B$- or $B^\cop$-Hopf module structure on $\lfr BV$
($\,\rfr VB\,$).\end{rem3}

\begin{rem4}[Sweedler \cite{Swe}] In the case of a Hopf algebra $H$ any left
(or right) Hopf $H$-module is left (resp. right) free, i.e. it has
the form $\lfr HV$ (resp. $\rfr VH$), with $V$ being (not
necessarily a finite dimensional) vector space of left (resp.
right) invariant elements.\end{rem4}

A left Hopf $B$-module which is at the same time a $B$-bimodule
satisfying $$ \cl{\lfr BV}(x.a)=\cl{\lfr BV}(x)\co(a)\eqno (20)$$
is called a left {\it covariant bimodule} \cite{Wor}. The right
$B$-module structure on $\lfr BV$ can be used to generate, via the
projection map (18), a right $B$-module structure on the vector
space $V$:
$$ \rho(a)v\doteq \ep_V((1\ot v).a) \ .\eqno (21) $$ Due to (18)
and (20) one obtains the following relationship $$(a\ot
v).b=ab_{(1)}\ot\rho(b_{(2)})v\eqno (22)$$ between right module
structures on $V$ and on $\lfr BV$. This means that the converse
statement is also true: any right $B$-module structure $\rho$ on
$V$ generates a left $B$-covariant bimodule structure on a left
free Hopf $B$-module $\lfr BV$. Similarly, the formula $$(a\ot
v).b=ab_{(2)}\ot\rho(b_{(1)})v\eqno (23)$$ induces a left
covariant $B^\cop$-bimodule structure on $\rfr VB$. In other
words, there is a bijective correspondence between right module
structures on $V$ and left covariant either $B$ or $B^\cop$ -
bimodule structures on $\lfr BV$. For a free right $B$-covariant
bimodule $\rfr VB$ one gets instead $$a.(v\ot b)=\lam(a_{(1)})v\ot
a_{(2)}b\eqno (24)$$ where $\lam$ denotes the left $B$-module
structure induced on $V$.

The following version of Lemma 4.1 is essentially due to
Woronowicz \cite{Wor}.
\begin{prop3}
Let $V$ be a finite-dimensional vector space and $B$ be a
bialgebra. Then  the following are equivalent:\begin{description}
\item[i)] A left $B$-module structure $\lam:B\rar\nd V.$
\item[ii)] A left $B$-module structure on a right free Hopf $B$-module
$\rfr VB$ such that it becomes a right (free) $B$-covariant
bimodule. Moreover, the commutation rule (16) takes the form
$\Lam^i_k (a)= \lam^i_k(a_{(1)})\,a_{(2)}$. Conversely,
$\lam^i_k=\ep\c\Lam^i_k$.
\item[iii)] A left $B$ (=$B^\cop$)-module structure on a right free Hopf
$B^\cop$-module $\rfr VB$ such that it becomes a right (free)
$B^\cop$-covariant bimodule. In this case, the commutation rule
takes the form $(\Lam^\cop)^i_k (a)= \lam^i_k(a_{(2)})\,a_{(1)}$
with $\lam^i_k=\ep\c(\Lam^\cop)^i_k$.
\end{description}\end{prop3}
\begin{proof} $iii)$ is a $B^\cop$ version of $ii)$ (22) and (23).
Taking into account (15) and (24) one calculates $$a.(e_k\ot 1)=
e_i\ot\Lam^i_k (a)= \lam^i_k(a_{(1)})e_i\ot a_{(2)}= e_i\ot
\lam^i_k(a_{(1)})\,a_{(2)}\ .$$ Hence $\Lam^i_k(a)=
\lam^i_k(a_{(1)})\,a_{(2)}$. Applying now $\ep$ to the both sides,
gives
$\ep\c\Lam^i_k(a)=\lam^i_k(a_{(1)}\,\ep(a_{(2)}))=\lam^i_k(a)$.
\end{proof}

\begin{rem5} Observe that the left-left $\yd$ condition (8) can be now
rewritten as $$\Lam^i_m (a)L^m_k = L^i_m(\Lam^\cop)^m_k (a) \ .$$
\end{rem5}

A similar statement holds true for left (free) covariant
bimodules. A $B$-bimodule which is at the same time a
$B$-bicomodule satisfying left and right Hopf $B$-module
conditions together with the left and right covariance  condition
is called a Hopf $B$-bimodule or, in the terminology of
Woronowicz, {\it bicovariant} bimodule.

Assume that $M\doteq\, \lfr BV$ is a left free bicovariant
bimodule. In this case, apart of the right multiplication (20) one
has at our disposal a right comultiplication $\ccr{M}:  B\ot V\rar
B\ot V\ot B$ such that
$$\,_M\!\co(a.x.b)\,=\,\co(a)\,_M\!\!\co(x)\co(b) \eqno (25)$$ and
the following bicomodule property $$(\id\ot\,\ccr M)\c\cl M= (\cl
M\ot\id)\c\,\ccr M \ . \eqno (26)$$
Here, $\cl M(a\ot v)=  a_{(1)}\ot a_{(2)}\ot v$ denotes a free
left comultiplication in $M$.
In this case, the right comultiplication $\ccr M$ in $M$ is
induced from a right comultiplication $\ccr V$ in $V$ according to
the formula $$_M\!\co\,(1\ot v)=1\ot\, _V\!\co (v)\eqno (27)$$ The
structure theorem (Drinfeld, 1986; Woronowicz, 1989; Yetter, 1991)
asserts that the vector space $V$, equipped with the right
multiplication (21) and the right comultiplication (27), inherits
a right crossed $B$-module structure. The inverse statement is
also true: a left (or right) crossed $B$-module structure on $V$
generates a left (resp. right) free Hopf $B$-bimodule structure on
$\rfr VB$ (resp. $\lfr BV$).

\begin{rem6} Due to Sweedler \cite{Swe} theorem, any bicovariant bimodule
$M$ over a Hopf algebra $H$ is free, i.e. it can be represented as
$\lfr HU$ or $\rfr WH$, where $U\equiv (U,\,
_U\!m,\,_U\!\co)\in\yd^H_H$ (resp. $W\equiv (W, m_W,\, \co_W
)\in\,_H^H\yd$) denotes a crossed bimodule of left (resp. right)
invariant elements in $M,$ Remark 4.3. If the antipode $S$ of $H$
is bijective then the following holds:\ $U\cong W$,
$_W\!m\,=(S^{-1})^\star(m_U)$ and $_W\!\co\,=S_\star\,(_U\!\co)$
(cf. Proposition 3.3).\end{rem6}

\section{BIMODULES DUAL TO FREE HOPF MODULES}
For an arbitrary left  $A$-module $M$ one can introduce its
$A$-dual: a right $A$-module $\ld M\doteq \mod (M, A)$ -  as a
collection of all left module maps from $M$ into $A$ (Bourbaki
1989). The evaluation map gives a canonical $A$-bilinear pairing
$$\forall\,x\in M\;\text{and}\;X\in\ \ld M,\quad\ll a.x,
X.b\gg\,\doteq\, a.X(x).b\in A.\eqno (28)$$

Similarly, for a right $A$-module $N$: $\rd N\doteq \mod(N, A)$, a
collection of all right module maps, seen as a left $A$-module, is
called an $A$-dual of $N$. In this case we shall write a canonical
pairing $$\forall\,y\in N\;\text{and}\;Y\in\rd N,\quad\ll a.Y,
y.b\gg\,\doteq\, a.Y(y).b\in A.\eqno (29)$$

For free finitely generated modules one can repeat the dual basis
construction.

\begin{lem2}[Bourbaki \cite{Bou}] Let $V$ be a finite-dimensional vector space.
$A$-dual module to the left (resp. right) free module $\lfr AV$
(resp. $\rfr VA$) can be represented as $\rfr {\tilde V}A$ (resp.
$\lfr A{\tilde V}$), ie.:$$\ld{(\lfr AV)}\,=\,\rfr {\tilde V}A\ \
\ \hbox{and}\ \ \ \rd{(\rfr VA)}\,=\,\ \lfr A{\tilde V},$$
$$\rd{(\ld{(\lfr AV)})}\,=\,\lfr AV\ \ \ \hbox{and}\ \ \
\ld{(\rd{(\rfr VA)})}\,=\,\rfr VA .$$
\end{lem2}

The canonical $A$-bilinear pairing $\ll\ ,\ \gg\,:\ \lfr A{\tilde
V}\ot\,\rfr VA\rar A$ between elements $a\ot\al\in\ \lfr A{\tilde
V}$ and $v\ot b\in\, \rfr VA$ can be rewritten, in this case, by
means of the usual $\bk$-bilinear pairing $\langle\ ,\ \rangle\, :
\tilde V\ot V\rar\bk$ $$ \ll a\ot\al, v\ot b\gg\, =\, ab\langle
\al, v\rangle\,\doteq\, ab\,\al(v) .\eqno (30)$$

Assume now that $B$ is a bialgebra and we have done a right
coaction $\ccr V: V\rar\,V\ot B$. The image of $\ccr V$ belongs to
a right free $B$-module $\rfr VB$. On the other hand, the
transpose left action $\widetilde{\ccr V}:\tilde V\rar\,\lfr
B{\tilde V}$ takes its values in $\lfr B{\tilde V}$\  --\ $B$-dual
to $\rfr VB$. This suggest a possibility for comparison between
both pairings:

\begin{lem3} Let $B$ be a Hopf algebra with antipode $S$.
Assume further that a finite dimensional vector space $V$ is a
right $B$-comodule with coaction $\ccr V$. Then for any $v\in
V$and $\al\in \tilde V$ one has $$1_B\langle \al, v\rangle\, =\,
\ll\widetilde{\ccr V}(\al), S_\star(\ccr V^\cop)(v)\gg \ .\eqno
(31)$$
\end{lem3}
\begin{proof} It is enough to check (31) on basis vectors: $$
\ll\widetilde{\ccr V}(e^k), S_\star(\ccr V^\cop)(e_i)\gg= \ll
R^k_j\ot e^j,\,e_m\ot S(R^m_i)\gg=R^k_jS(R^j_i)=1_B\de^k_i$$ due
to (4), (7), (30) and using the properties of the antipode.
\end{proof}

 If $M$ is $A$-bimodule then $\ld M\doteq\hom_{(A,-)}(M, A)$
can be equipped in a canonical bimodule structure \cite{Bor1,Bor2}
by $$\ll x, a.X.b\gg\,\doteq\, X(x.a).b\,=\,\ll x.a, X\gg .b \eqno
(32)$$ We call the bimodule $\ld M$ a left $A$-dual of $M$
Similarly, one can define a right dual $\rd
M\doteq\hom_{(-,A)}(M,\, A)$.

Let $\rfr VA$ be a right free bimodule with a left module
structure given by the commutation rule (15). Then its right
$A$-dual $_A\tilde V$ is a left free bimodule with the transpose
commutation rule, $$ \ll 1\ot e^k,\,a.(e_i\ot 1)\gg\, =\,\ll 1\ot
e^k,\, e_m\ot \Lam_i^m(a)\gg\,=$$ $$=\,\langle e^k, e_m
\rangle\,\Lam^m_i(a)\, =\,\Lam_i^k(a),$$ $$\ll 1\ot
e^k,\,a.(e_i\ot 1)\gg\,=\ll(1\ot e^k).a,\,e_i\ot
1\gg\,=\ll\Phi^k_m(a)\ot e^m,\,e_i\ot 1\gg\,=\,\Phi^k_i(a)$$
Therefore, $\Phi=\tilde\Lam$. Assume further that $A$ is a
bialgebra, hence $\Lam$ has the form (24)
$$\Lam(a)=\lam(a_{(1)})\,a_{(2)} \eqno (33)$$ for some
representation $\lam$ of $A$ in $V$. This means that $\rfr VA$ is
a right $A$-covariant Hopf bimodule generated by $(V, \lam)$. Thus
the transpose commutation rule $\Phi(a)= \tilde\lam
(a_{(1)})\,a_{(2)}$ makes $\lfr A{\tilde V}$ a left
$A^\cop$-covariant Hopf bimodule (cf. (23) and Appendix).

These prove our main result:
\begin{thm} Assume $V$ is a finite-dimensional vector space and $A$ is an
algebra.
\item[i)] Let $\rfr VA$ be a right free $A$-bimodule whose left
module structure is given by a commutation rule $\Lam: A\rar
A\ot\nd V$. Then its right $A$-dual $\rd {(\rfr VA)}=\,\lfr
A{\tilde V}$ is a left free $A$-bimodule with the transpose
(right) commutation rule $\Phi=\tilde\Lam$.
\item[ii)] Assume further that $A$ is a bialgebra and $\rfr VA$
is a right $A$-covariant bimodule generated by representation
$\lam: A\rar\nd V$ (33). Then its right $A$-dual $\lfr A{\tilde
V}$ is a left (free) $A^\cop$-covariant bimodule generated by
$(\tilde V,\,\tilde\lam)$.
\end{thm}
The above theorem suggests that dualizing a bicovariant bimodule
over a bialgebra $B$  one obtains, in general, a $B^\cop$-
covariant bimodule which is not necessarily bicovariant (unless
$B\cong B^\cop$). However, in the special case of quantum groups
we get

\begin{cor2} Let $B$ be a bialgebra with a bijective antipode (quantum
group). Let $\lfr BV$ be a bicovariant $B$-bimodule  generated by
a right-right $\yd$ module $(V,\,_Vm,\,_V\!\co)\in\,\yd^B_B$. Then
its left $B$-dual  $\rfr {\tilde V}B$ is a right (free)
$B^\cop$-covariant bimodule generated by $(\tilde
V,\,\widetilde{_V\!m})$.  Due to Corollary 3.4 and Proposition
3.5, it can be also equipped,  with a bicovariant
$B^\cop$-bimodule structure generated by $(\tilde
V,\,\widetilde{_V\!m},\,(S^{-1})_\star\,(\widetilde{_V\!\co}^\cop)
\in\,^{B^\cop}_{B^\cop}\yd$. Thus the identity (31) is satisfied.
\end{cor2}

This observation can be relevant for adapting  the vector fields
formalism \cite{Bor1,Bor2} to the case of differential calculus on
quantum groups \cite{Wor}. It also corrects the claim formulated
in Aschieri \& Schupp \cite{AS,Asc} that general vector fields for
a bicovariant differential calculus form a bicovariant bimodule
over the same Hopf algebra.

Consider a right-covariant differential calculus $d:B\rar\Gamma$
over a bialgebra $B$ with values in a free right $B$-covariant
bimodule $\Gamma$ of one-forms. Thus $\Gamma\cong\,\rfr VB$ and
the left $B$-module of right invariant elements $V\equiv (V,
\lam)$ generates the left $B$-module structure of $\Gamma$. The
right-covariance property writes \cite{Wor}: $$_\Gamma\!\co\c d =
(d\ot\id)\c\co\ ,\eqno (34)$$ i.e. $d:\Gamma\rar B$ is a right
comodule map. The right dual $\Gamma^\dagger\cong\,\lfr B\tilde
V$, where the generating space $(\tilde V, \tilde{\lam})$ is a
right $B$-module, is a left $B^\cop$-covariant bimodule. The
generalized Cartan formula $$X^\pt (f)\doteq\ll X, df\gg \eqno
(35)$$ allows us to associate with any element
$X\in\Gamma^\dagger$ the corresponding $\bk$-linear endomorphism
$X^\pt\in\nd\!_{\bk} B$. This gives the action
$\pt:\Gamma^\dagger\rar\nd\!_{\bk} B$ which satisfies the axioms
of a right Cartan pair \cite{Bor1,Bor2}.
\begin{prop4}
With the assumptions as above, for any element $\al\in\tilde V$
the corresponding endomorphism $(1\ot\al)^\pt:B\rar B$ is a right
comodule map. More exactly  one has
$$\co\c(1\ot\al)^\pt=((1\ot\al)^\pt\ot\id)\c\co\ .\eqno (36)$$
\end{prop4}
\begin{proof} Let $\{e_i\}$ be any basis in $V$ and $\{e^i\}$ the
dual basis in $\tilde V$. With respect to the given basis one can
define the generalize derivations $\pt^i\in\nd\!_{\bk}B$ as
$$\pt^i\doteq (1\ot e^i)^\pt\ . $$ Thus $\pt^i(fg) =
\pt^i(f)g+\lam^i_k(f_{(1)})f_{(2)}\pt^kg$. Substituting
$df=e_i\ot\pt^if$ into equation (34) and comparing the
coefficients in front of the same basis vectors we conclude
$$\pt^if_{(1)}\ot f_{(2)}=(\pt^if)_{(1)}\ot (\pt^if)_{(2)} \ \ \
\text{for}\ \ \forall\, i\ ,$$ which is equivalent to (36).
\end{proof}

Let now $B$ be a Hopf algebra with bijective antipode. Consider a
Woronowicz bicovariant differential calculus $d: B\rar\Gamma$,
where $\Gamma$ is $B$-bicovariant bimodule of one-forms. Thus
$\Gamma\cong\,\rfr VB$ and $(V,\,m_V,\,\co_V)\in\,^B_B\yd$ denotes
a crossed bimodule of the right invariant elements of $\Gamma,$
Remark 4.6. Dualizing the bimodule of one forms, one obtains a
$B^\cop$-bicovariant bimodule $(\Gamma)^\dagger\cong\,\lfr
B{\tilde V}$ of vector fields where now, by Corollary 3.4,
$(\tilde V,\,\widetilde{m_V},\,(S^{-1})_\star\,
(\widetilde{\co_V}^\cop) \in\,\yd^{B^\cop}_{B^\cop}$ denotes the
crossed module of the Woronowicz left invariant vector fields.
More exactly, the Woronowicz's vector fields can be identify with
functionals on $B$ spanned by $$\chi^i\doteq \ep\c\pt^i\ .$$ The
quantum Lie bracket structure is induced by a convolution product
and the Yang-Baxter operator on $\tilde V$.

\section*{Acknowledgements} The authors  thank Micho
\mbox{{\Dj}ur{\dj}evich}, Piotr Hajac and W{\l}adys{\l}aw Marcinek
for discussions.

\section*{Appendix}
We give here an alternative, i.e. by direct calculations, proof of
Theorem 5.3 \textit{ii)}. For this aim, it suffices to check that
the right multiplication (cf. (23-24))$$(1\ot\al).a\doteq
a_{(2)}\ot \tilde \lam(a_{(1)})\al \eqno{(A1)} $$ and the left
comultiplication (see Remark 4.2) $$\cl{\lfr A\tilde
V}(b\ot\al)\doteq b_{(2)}\ot b_{(1)}\ot\al \eqno{(A2)}$$ in $\lfr
A\tilde V$ are related by the left $A^\cop$-covariant condition
(cf. (20)) $$\cl{\lfr A\tilde V}((1\ot\al).a)=\cl{\lfr A\tilde
V}(1\ot\al) \co^\cop(a).\eqno{(A3)}$$ Calculation of the left-hand
side gives $$\cl{\lfr A\tilde V}(a_{(2)}\ot \tilde
\lam(a_{(1)})\al)= (a_{(2)})_{(2)}\ot (a_{(2)})_{(1)}\ot \tilde
\lam(a_{(1)})\al \ .$$ From the other hand similar calculations
for the right-hand side yield $$
a_{(2)}\ot(1\ot\al).a_{(1)}=a_{(2)}\ot (a_{(1)})_{(2)}\ot \tilde
\lam((a_{(1)})_{(1)})\al \ .$$ The proof is finished since
$$(a_{(2)})_{(2)}\ot (a_{(2)})_{(1)}\ot a_{(1)}= a_{(2)}\ot
(a_{(1)})_{(2)}\ot (a_{(1)})_{(1)}$$ due to the coassociativity
property;
$(\co^\cop\ot\id)\c\co^\cop=(\id\ot\co^\cop)\c\co^\cop$.\hfill$\Box$

\end{document}